\documentclass[10pt,oneside,a4paper]{article}

\usepackage{amscd,amssymb,amsmath,amsthm,mathrsfs,dsfont}
\usepackage[a4paper]{geometry}
\usepackage{color}
\usepackage{tikz} 
\usetikzlibrary{positioning}
\usepackage[shortlabels]{enumitem}
\usepackage{comment}
\usepackage{caption}
\usepackage{array}
\usepackage{dirtytalk}
\usepackage{subcaption}
\usepackage[thinc]{esdiff}

\usepackage{url}
\usepackage[hyperindex,breaklinks]{hyperref}
\hypersetup{
    colorlinks=true,
    linkcolor=blue,
    filecolor=black,      
    urlcolor=blue,
    citecolor=blue,
}
  \hypersetup{
           breaklinks=true,   
           colorlinks=true,   
        }
\usepackage[hyphenbreaks]{breakurl}
\def\Z{\mathbb{Z}}
\def\R{\mathbb{R}}
\def\P{\mathbb{P}}

\numberwithin{equation}{section}
\usepackage{fancyhdr}
\usepackage[backend=bibtex,style=ieee,giveninits=true,sorting=anyt]{biblatex}
\newcommand\shorttitle{The Ising model: Highlights and perspectives}

\begin{document}

\hypersetup{linkcolor=black}
\hypersetup{urlcolor=black}
\title{The Ising model: Highlights and perspectives}
\author{Christof K\"ulske\footnotemark[1] 
}
\date{\today}

\maketitle
\begin{abstract} We give a short non-technical introduction to the Ising model, and review some successes as well as  
challenges which have emerged from its study in probability and mathematical physics. This includes the infinite-volume theory of phase transitions, and ideas like scaling, renormalization group, universality, SLE, and random symmetry breaking in disordered systems and networks. \medskip

This note is based on a talk given on 15 August 2024, as part of the Ising lecture during the 11th Bernoulli-IMS world congress, Bochum. 
\end{abstract}
\textbf{AMS 2000 subject classification:} 60K35, 82B20, 82B26\vspace{0.3cm}

{\bf Key words:}  Ising model,  spin models, Gibbs measures, phase transitions, renormalization group, SLE, disordered systems, extremal decomposition.    
\footnotetext[1]{Ruhr-University Bochum, Germany}
\footnotetext[0]{E-Mail: \href{mailto:Christof.Kuelske@ruhr-uni-bochum.de}
{christof.kuelske@rub.de}; ORCID iD: \href{https://orcid.org/0000-0001-9975-8329}{0000-0001-9975-8329}}
\footnotetext[0]{\url{https://math.ruhr-uni-bochum.de/fakultaet/arbeitsbereiche/stochastik/gruppe-kuelske/mitarbeiter/christof-kuelske/}}
\tableofcontents
\hypersetup{linkcolor=blue}

\newpage

\section{A theory for ferromagnetic ordering?}

Ernst Ising  (1900 Cologne - 1998 Peoria, Illinois) was born into a Jewish family from Cologne 
who moved to Bochum, an industrial city in Ruhr area in the western part of Germany. Ernst Ising grew 
up in Bochum, obtained his Abitur (high-school diploma) in Bochum and went on to study physics 
in G\"ottingen, Bonn and Hamburg. 
In Hamburg he followed Ph.D. studies with Wilhelm Lenz (1888 Frankfurt - 1957 Hamburg) as  a 
Ph.D. advisor. Ising later escaped the prosecution during the Nazi period 
via Luxemburg and settled in the US.  

The problem of Ising's Ph.D. thesis was  the following: 

\textit{Can the formalism of statistical mechanics, 
following Boltzmann and Gibbs, explain a phenomenon like ferromagnetism?}  

\textit{Ferromagnetism} is the physical property of a  magnetic material (e.g. iron) 
that it can be magnetized by putting on an external magnetic field at low temperatures. If the magnetic field 
is then decreased to zero, the magnetization does not vanish, and one sees the so-called \textit{spontaneous 
magnetization}. If then temperature is slowly increased, the spontaneous magnetization decreases 
and vanishes sharply at a particular finite \textit{critical temperature} whose value depends on the material. 

\textit{Statistical mechanics} describes equilibrium states of large systems 
by means of weighted averages (instead of trying to solve 
time-dependent evolution equations for many-particle systems in detail). 
 The weights are depending on the individual states of all particles. They are  
given by exponentials of the 
negative total energy of the system, with a prefactor $\beta>0$ 
whose meaning is the inverse of the temperature. An excellent introductory book to the subject is \cite{FrVe17}.

\section{Definition of the Ising model} 

It was known at the time of a young Ising 
that metals have a spatially regular structure, with atoms sitting 
on regular lattices. 
Based on this, the following statistical mechanics model which was later called \textit{Ising model} for ferromagnetism 
was invented in 1920. The idea was to simplify the description of the 
local degrees of freedom (\textit{local magnetic moments}) in an extreme way, and further 
assume only an interaction between nearest neighbors on the lattice of such magnetic moments, sitting on the 
sites of a lattice $\Z^d$. The spatial dimension $d$ should in  physical applications mostly take the values 
$3$ or $2$, but mathematically $d$ may also be an arbitrary integer.  \\

In the Ising model, the degrees of freedom of the (finite-volume) model are given by  \textit{magnetic moments (spins)} 
$\sigma_i\in \{1,-1\}$ which are sitting at sites $i\in \Lambda \subset\Z^d$, where $\Lambda$ is 
a very large but finite subset of the lattice. 
 In today's probabilistic language they are random variables on a common probability space, which 
 are distributed according to the following \textit{Gibbs distribution}. 
This is the probability measure depending on inverse temperature $\beta>0$, and external field $h\in \R$ 
in which the probability to see a configuration $\sigma_{\Lambda}=(\sigma_i)_{i\in \Lambda}$ 
is given by $$\frac{1}{Z_{\Lambda}}\exp\big(\beta \sum_{i \sim j} \sigma_i \sigma_j + h\sum_{i}\sigma_i \big), $$  
where the first sum is over nearest neighbors $i \sim j$ in the finite volume $\Lambda$.  
The negative exponent is (up to constants) the energy of a configuration.  Spin configurations 
which for most edges $i \sim j$ align 
are preferred under this measure, as the product of two spins $\sigma_i \sigma_j$ is maximal iff they are equal.  

Next to the lattice $\Z^d$, other choices of base spaces which are very interesting 
today  are \textit{networks}. Then the model is used also with other interpretations of the "local magnetic moment" $\sigma_i$ in mind,  
e.g. as \textit{opinion of a person $i$ in a network.}   

The above exponential formula for the Gibbs measure in finite volume 
has free boundary conditions. In principle there may (and often should) be terms 
connecting the sites which sit on an inside layer of  the finite volume to a fixed boundary condition. 
These couple to the spins inside only via its values on the outside boundary layer of the volume.  
Such a boundary condition could for instance 
be chosen to be the maximal boundary condition, in which all spins on the 
sites on the outer boundary are equal to one.  
Now, does this model show a \textit{phase transition} and \textit{long-range order} (in which 
spins tend to align over long distances)? \medskip

\noindent $\bullet$ Is there a spontaneous magnetization at large $\beta$? \\ 
$\bullet$ Does the model have a magnetic  phase transition? \\ 
$\bullet$ What is actually a phase transition in a precise mathematical sense? \\ 
$\bullet$ What can be said about the behavior of other physical quantities?\medskip 

Other physical quantities of interest include the specific heat (the amount of energy needed to 
put into the system in order to increase temperature) and the susceptibility (the increase of the total 
magnetic moment of the material under increase of the external magnetic field). \\ 

Ernst Ising was the first person on earth to begin to tackle such questions. 
He defended his Ph.D. thesis in the year of 1924 in Hamburg, in which 
he showed by explicit combinatorial computations  
that the model has no phase transition in spatial dimension $d=1$. 
The result was published in \textit{Zeitschrift f\"ur Physik}  (since 1997 European Physical Journal) 
in the year of 1925 \cite{Is1925}.
He also conjectured in this paper that the absence of a phase transition in the short range model 
would probably be true also in higher spatial 
dimensions, and long-range interactions would probably be necessary to create a spontaneous magnetization 
in which spins over long distances stay aligned. 
Luckily he was too pessimistic about the short-range Ising model. 
Independently we note that the study of long-range models (in the form of Dyson models with power-law decay of the strength 
of the  interactions between spins at sites $i$ and $j$, as a function of the distance $d(i,j)$) 
turned out to be a beautiful story in itself which developed much later,  
where ferromagnetic ordering may even appear in dimension $d=1$, 
see \cite{Dy69, FrSp82,AiChChNe88}.

\section{Phase transitions and long-range order in lattice dimensions \texorpdfstring{$d\geq 2$}{Lg}}

{\bf Boundary condition persistence via Peierls contours and Peierls argument. } Rudolph Peierls (1907 Berlin - 1995 Oxford) was a German-Jewish physicist who grew up in Berlin, and 
then moved to England.  
Already in 1936 he came up with his famous Peierls argument \cite{Gr64,FrVe17}
which gave a proof that, contrary to 
the conjecture of Ising,  there actually was long-range 
order in the two-dimensional short-range Ising model, in zero external field, at sufficiently low temperatures. 

His statement was about the finite-volume Gibbs measure in large boxes, with an all-plus boundary condition. 
He proved that the probability to see a spin value $\sigma_v=-1$ at a site in the middle of the box, 
under this measure, is bounded above by $e^{- c\beta}$ for $\beta$ large enough, with a bound which holds 
uniformly even for arbitrarily large boxes. In other words, the influence of the boundary condition 
persists over arbitrarily large distances, which indicates that there is \textit{long range order} 
in the Ising model.  

His proof idea is based on the observation that in each spin configuration in the box which has $\sigma_v=-1$, 
there must be at least one contour surrounding it, which separates plus clusters of spins from minus 
clusters of spins.  We call such contours today \textit{Peierls contours}. It is then easy to bound the 
probability of a fixed contour of length $l$ by $e^{- 2 \beta l}$. Together with a bound 
of the number of length-$l$ contours around a given site of the form $ C l 3^l$, this proves the claim when 
summing over the possible lengths $l$ of contours.   

Versions of Peierls contours and versions of the Peierls argument 
are used today in many places different from the original nearest neighbor Ising 
model. This includes percolation theory \cite{Gr99, BoRi06} (which studies the appearence of large connected clusters in spatial
random models) and 
spin models with long-range interactions. 

Contours are also key objects in \textit{Pirogov-Sinai theory} \cite{Si82,Za84, EnFeSo93, HePeRe20}, with contributions 
from Yakov Sinai (born 1935 in Moscow, Abel prize 2014).  
This is a general theory 
to analyze \textit{phase diagrams} of low temperature lattice spin models which are 
 non-symmetric generalizations of the Ising model to situations where the 
 spins may take many spin values. 
A phase diagram is a decomposition of the parameter space (of temperature and finitely many external fields) 
into regions in which the structure and number of macroscopic states (see Subsection \ref{sec:InfiniteSystems}) does not change.   
Pirogov-Sinai theory is based on suitably adapted complex series expansion methods (cluster expansions)
for logarithms of complex polynomials in many variables (given by possibly complex-valued contour weights), 
whose coefficients are governed by combinatorial rules. 

The idea of Peierls to study boundary curves (interfaces) separating plusses from minusses in the Ising model 
turned out to be very fruitful in a number of different places, and in particular 
much later in history as a starting point to analyze limiting behavior of random curves 
which led to stochastic Loewner evolution (SLE), see below Section \ref{sec:SLE}.
\medskip 

{\bf Non-analyticity of free energy as a function of temperature via exact Onsager solution. } 
Lars Onsager (1903 Oslo - 1976 Florida) was a theoretical chemist who  
received the Nobel prize in chemistry in the year of 1968 for 
his contributions to non-equilibrium statistical mechanics  (where one studies e.g. systems 
at not just one temperature, but in contact with reservoirs with two different temperatures) \cite{On31}. 

A very highly valued achievement in the realm of thermodynamic equilibrium was the so-called Onsager-solution of 
the two-dimensional Ising model, which was published in Phys. Rev. 1944 \cite{On44}, see also \cite{CoWu73}. \textit{The solution of the model}  
in this context refers to the computation of  the \textit{free energy} of the Ising model in zero external field $h=0$ in lattice dimension $d=2$. Onsager  was able to give the explicit (and correct) formula 
\begin{equation*}\begin{split}\label{eq:DL:joint}
&\lim_{\Lambda \uparrow \Z^2}\frac{1}{|\Lambda|}\log Z_{\Lambda}(\beta)=\log 2 \cr
&+ \frac{1}{8 \pi^2}\int_{0}^{2\pi}\int_{0}^{2\pi}\log \Bigl( 
\cosh(2 \beta)^2-\sinh(2 \beta)\bigl( \cos \theta_1 + \cos \theta_2 \bigr) \Bigr)
d\theta_1 d\theta_2
 \end{split}
\end{equation*}
where one take the  volume limit $\Lambda \uparrow \Z^2$ (for instance) along sequences of boxes with increasing sidelength which are all 
centered at the origin. 
To arrive at this expression 
Onsager used graphical representations of the partition functions $Z_{\Lambda}(\beta)$ of the model, 
combinatorics, algebra, determinants (respectively Pfaffians which are square-roots of determinants). 
Modern versions of the proof proceed via so-called fermionic integration (defined in terms 
of projections of Grassmann algebras) to handle combinatorics and to derive this formula. 

Once this formula is obtained, it is a simple matter to see that there is singular behavior as a function 
of the inverse temperature $\beta$ only at that critical value for which the argument of the $\log$ touches 
zero when both cosines equal $1$ which amounts to  $2\beta_c=\log(1+\sqrt{2})$. 
Probabilistically,  the free energy is the log-moment generating function of 
the (negative) energy $\sum_{i \sim j} \sigma_i \sigma_j $ per lattice site. Taking derivatives w.r.t. temperature
one obtains the famous \textit{logarithmic 
divergence of the specific heat} of leading order  
$-\log| \beta-\beta_c|$ when $\beta$ approaches the critical inverse temperature. 

By a further ingenious piece of work, which remained unpublished for a long time,  
Onsager and Kaufmann showed that there is power law critical behavior for magnetisation \cite{Ba11}. 
Via intricate computations 
they found that the spontaneous magnetization per site behaves like 
 $$m(\beta)=\lim_{n}\langle \sigma_0 \rangle_{\Lambda_n,+} \sim (\beta - \beta_c)^{\frac{1}{8}}$$ 
 as $\beta \downarrow \beta_c$. 
Here  the brackets denote the expectation of the spin variable at the origin in the measure with plus boundary 
 conditions in a box $\Lambda_n$ which is centered at the origin, as the sidelength tends to infinity. 
 
 The planar Ising model is not the only statistical model which admits a solution formula in terms of some determinant \cite{Ba82,KeOk06},  but solvability in this sense seems to be possible only in spatial dimension $d=2$.
 

\section{Critical behavior and universality: \texorpdfstring{\\} a physical perspective} 

{\bf Critical exponents. }  Can the Ising model really be relevant for the explanation of the behavior of the physics 
of a real material, 
with all the extreme simplifications going into the formulation of the model? 
Surprisingly, yes. Many materials behave like prototypical models,  
 e.g. like the Ising model. 
Namely, physical experiments with real magnets performed around the 1970s
support the following so-called \textit{universality conjecture}:   \medskip

\textit{Critical behavior and in particular critical exponents for systems with short-range interactions 
depend only on few general properties like 
spatial dimension $d$ and dimension of the local spin space.}   \medskip

There are relevant models in which the Ising local state space (or fibre) $\{-1,1\}$ is replaced 
by a copy of $\R^k$, or a manifold like a sphere,  
in which cases one sees critical behavior different from that of the Ising universality class. 
Here \textit{critical behavior} refers to the leading asymptotic behavior of physical quantities like the magnetization in a 
neighborhood of the critical point $\beta=\beta_c$. Many substances show 
power law behavior of a form we just saw for the special case of the planar Ising model, where the 
\textit{magnetization critical exponent} takes the exact value $\frac{1}{8}$.  Today it is known that the magnetization 
in the Ising model behaves like 
$$
m(\beta) \sim \begin{cases}
			(\beta - \beta_c)^{0.3264...}  & d=3\\
            (\beta - \beta_c)^{\frac{1}{2}}  & d\geq 4 		 \end{cases}
$$
where the asymptotics in the last line ignores logarithmic corrections in $d=4$. 
The exponent for $d=3$ is supported by physical experiments, numerical simulations, and 
non-rigorous approximation schemes (which nevertheless seem to work very well). 
The critical exponents for large spatial dimensions $d$ coincide with those of \textit{mean-field models} 
which are models on the complete graph with vertex set $\{1,2,\dots,n\}$ where all possible 
pairs of spins contribute an interaction $\beta \sigma_i \sigma_j$ to the Boltzmann weight. 

But this is hard to show. The common but technically involved approach is based on the lace expansion, a diagrammatic expansions
which is able to prove mean-field 
behavior in sufficiently spread out models, or nearest neighbor models in sufficiently high dimensions \cite{BoHoKo18,Sl06}. 
A related angle is to prove Gaussianity 
of the scaling limit as in \cite{AiCo21}.

\medskip

How to explain such universality of the lattice Ising model with its dimensional dependence even heuristically? 
Universality is well-known in the form of classical limit theorems in probability. It is best known in the form of 
the \textit{central limit theorem} \cite{Kl14} which holds 
for standardized sums of i.i.d. random variables when their number tends to infinity, 
and correspondingly distributional convergence of random walk paths to Brownian motion. 

{\bf Renormalization group: flow heuristics, non-localities and singularities, successes, problems. }  This suggests to look for non-standard limit theorems for 
the lattice-indexed dependent Ising spin-variables, under a suitable scaling (zooming out). 
A non-rigorous but brilliant breakthrough in this regard was made by Kenneth Wilson
(1936 Massachusetts - 2013 Maine) who received the Nobel prize in theoretical physics for precisely this work 
in the year of 1982. He outlined how 
universality and power laws could be explained by (non-rigorous) renormalization group theory \cite{Wi75}. 
Renormalization comes in various meanings in theoretical physics.  Historically even earlier forms of renormalization 
derive from the study of series expansions needed in particle physics \cite{AmMa05,Zi89}. 
In theoretical physics, renormalization group (RG) is associated with a recipe for handling certain (typically non-convergent) series, which consists of 
formal book-keeping rules in the form of diagrams which depict the countable zoo of terms appearing
in the corresponding expansion, and adhoc rules to remove infinite terms. Many of these steps are hard to justify mathematically.   

Wilson-renormalization in statistical mechanics now comes 
with the idea that carrying out scaling transformations (coarse-graining away fine scale details of configurations, and zooming out) should  correspond to a so-called \textit{renormalization group flow} 
in a space of finitely many parameters, here generically denoted by $p$,  
in some suitable space of interactions. 
This non-linear renormalization group flow $(t,p_0)\mapsto \Phi(t,p_0)$ 
depends on a real parameter $t \geq 0$ which describes how much we are zooming out, starting from parameter values $p_0$,  
and is expected to display fixed points with stable and unstable manifolds attached to it.  The asymptotics 
locally around these fixed points 
should correspond to the  type of critical behavior for a whole universality class and thus determine 
the critical exponents of the original system. 
In Wilson-theory \cite{Wi75} (to deal with the Ising model in $d=3$ for instance)  
this is accompanied with a number of further mathematically non-justifiable tricks for actual computations, 
one of which is expanding and 
taking formal limits in the dimensional parameter $d$, thereby treating it as if it were a continuous variable. \medskip 


Mathematical physics has not been able to justify Wilson-renormalization in its literal form. 
This poses many open challenges! 
To begin with one major problem, due to the so-called \textit{RG-pathologies}, RG transformations 
tend to create very long-range interactions, which can not simply be ignored \cite{EnFeSo93, Gr79}. Lattice systems tend 
to leave any reasonable space of interactions immediately under coarse-graining.

In order to fully understand this phenomenon, one needs to go to a proper infinite-system description, see Section \ref{sec:InfiniteSystems}. 
Indeed, non-localities and loss of regularity of the measures appear more generally also for many other types of 
transformations of lattice-indexed distributions of random variables, 
including stochstic time-evolutions. There is an underlying physical mechanism 
for this production of singularities, which is related to \textit{hidden phase transitions} in those 
parts of the systems which are integrated out. Depending on the concrete system and the transformation
a multitude of scenarios and forms of singularities may appear. 
Let us mention e.g. deterministic spatial projections \cite{Sc89}, 
coarse-graining transformations in local state space \cite{MaReShMo00}, applications of noisy kernels and stochastic time-evolutions \cite{EnFeHoRe02, KuLe07}, joint measures of spin variables and disorder variables \cite{KuLeRe04}
in \textit{quenched disordered systems} which 
we will discuss in the last section.

However, renormalization group and scaling, taken with a broad view, for instance on the level of contour and interfacial 
lines, 
have been guiding principles for beautiful fully rigorous mathematical physics and probability.  
We mention the study of low-temperature random systems \cite{BrKu88}, 
and the  description of two-dimensional critical system in terms of random curves with classes of distributions 
given by  $SLE_{\kappa}$, see below 
in Section \ref{sec:SLE}. There is 
application of Wilson renormalization group ideas to \textit{stochastic 
PDEs} \cite{Ku16}, and also 
relation to the theory of 
regularity structures \cite{Ha14} to deal with fluctuations and possible divergences
on a local scale 
for which Martin Hairer was awarded the Fields medal in 2014. 

\section{What is a phase transition? Infinite systems and Gibbs simplex} 
\label{sec:InfiniteSystems}

It is a fruitful idealization to study, instead of \textit{very large finite} systems, systems which are actually 
\textit{infinitely} large. The formalism we are using for this today goes back to 
Roland Dobrushin (1929 St Petersburg - 1995 Moscow) who obtained his Ph.D. with advisor Kolmogorov in 1955, 
Oscar Lanford  (1940 New York City -  2013 Switzerland), and David Ruelle 
(1935 born in Ghent, Belgium).   
The mathematical theory of infinite systems has many advantages, as phase transitions become sharp, 
and there is a natural way to describe phase coexistence for the same model parameters 
in terms of \textit{multiple Gibbs measures}. A canonical reference 
for large parts of the theory discussed in this section
is the clear 
and abstract  book by Hans-Otto Georgii \cite{Ge11}.

{\bf DLR-formalism. }  In the infinite volume, as Dobrushin, Lanford and Ruelle \cite{Ru99} realized, 
a notion of  \textit{consistency of the measure} has 
to be defined, as it is impossible to directly write expressions for Gibbs probabilities involving infinitely many sums, 
and quantities like the partition function $Z_{\Lambda}$ are only finite for finite $\Lambda$. 
This DLR-consistency equation has analogy to the much simpler situation for Markov chains. 
For Markov chains, one has only one determining equation for its 
invariant distribution, when the Markov transition matrix is given. The invariant distribution is then called stationary 
distribution, i.e. it is invariant under the associated stochastic 
time evolution. The more sophisticated analogue to a single  transition matrix, 
is in DLR theory given by a whole family of so-called specification kernels $\gamma_{\Lambda}$, 
where $\Lambda$ runs over the finite subvolumes of $\Z^d$. 
These kernels are precisely the finite-volume Gibbs measures we already encountered, 
but with a boundary condition $\omega$ 
outside of $\Lambda$. Thus, for the Ising model they take the form  
$$\gamma_{\Lambda}(\sigma_{\Lambda}| \omega_{\Lambda^c})
=\frac{1}{Z_{\Lambda}(\beta,h, \omega)}
\exp\Bigl(\beta \sum_{\substack{i \sim j\\ i,j \in \Lambda}} \sigma_i \sigma_j + h\sum_{i\in \Lambda}\sigma_i  + \beta \sum_{\substack{i \sim j\\ i\in \Lambda, j \in \Lambda^c}} \sigma_i \omega_j \Bigr) $$  
Here the fixed choice of boundary condition $\omega$ specifies a probability distribution for the variable $\sigma$ 
inside $\Lambda$. The above expression thus is a kernel, as it describes a probability measure w.r.t. $\sigma$ which 
depends measurably on the boundary condition $\omega$. 

Let now $\mu$ be a candidate infinite-volume state, i.e. a probability measure on $\{-1,1\}^{\Z^d}$, together with the standard product sigma-algebra on this space of infinite configurations.  
$\mu$ is called a \textit{Gibbs measure (in the infinite volume)} for a statistical mechanics model described by a family of 
\textit{specification kernels} $\gamma_{\Lambda}$, if and only if it satisfies the \textit{DLR consistency equation}
$$\mu \gamma_{\Lambda}=\mu$$ for all finite volumes $\Lambda$. 
The left hand side of the DLR equation describes the application of the kernel $\gamma_{\Lambda}$ 
to the measure $\mu$, 
by taking a boundary condition $\omega$ from $\mu(d\omega)$ and re-randomizing $\sigma_{\Lambda}$ accordingly 
inside the finite volume $\Lambda$.  
The DLR equation demands that this procedure should leave $\mu$ invariant as a measure, 
and this invariance under re-randomization should hold for all finite volumes $\Lambda$. 
The DLR formalism is different from the description of an infinite-volume measure via the Kolmogorov extension theory, 
as there may be multiple solutions $\mu$.  In general 
their finite-dimensional marginals can not in an obvious way directly read off from the DLR equation. \medskip

{\bf Gibbs simplex, extremal decomposition, pure states.} The DLR equations are linear in $\mu$, hence convex combinations $\alpha \mu_1+(1-\alpha)\mu_2$ 
of solutions $\mu_1,\mu_2$ are again solutions. 
Consequently the set of Gibbs measures is a convex set of probability measures in the infinite volume. 
It turns out that it is even a \textit{simplex}, which means that each Gibbs measure 
has a \textit{unique} decomposition $\mu=\int\nu \rho_{\mu}(d\nu)$ as an  integral over the \textit{extremal (non-decomposable) Gibbs measures}, with some 
decomposition measure $\rho_{\mu}(d\nu)$.  Integrals over sets of measures like this are measure theoretically 
well-defined by putting the natural sigma-algebra on the space of infinite-volume measures, the so-called evaluation sigma-algebra. 
This allows to apply the previous formula to any spin-observable. 
Extremal Gibbs measures are also called \textit{pure states} or \textit{pure phases}, 
and they are characterized by the fact that they are \textit{tail-trivial}, i.e. trivial on the tail-sigma algebra of spin-events. This is 
a generalization of the familiar Kolmogorov zero-one law for independent random variables. 
It means that a tail-event 
(that is an event in spin space which does not depend on any change of the values of finitely many spins) 
obtains only probability zero or one, in any pure state of a statistical mechanics model described 
by a specification.  An example of a tail-measure observable in the Ising model 
is the \textit{magnetization per site} 
$\limsup_{n}\frac{1}{|\Lambda_n|}\sum_{i \in \Lambda_n}\sigma_i $ along a sequence of finite volumes $\Lambda_n$ 
tending to infinity, or more generally empirical sums of some local observable of interest. The zero-one law implies that such a 
macroscopic quantity has a sharp non-fluctuating value in any pure state. Tail-triviality is mathematically equivalent
to spatial correlation decay of observables, in a suitable sense. 
For a clear exposition from a probabilistic point of view, 
see \cite{Ge11}, and also 
\cite{RaSe15, Si93}.

What is now the meaning of a pure state in infinite-volume statistical mechanics for physics?  The \textit{mathematical pure states} are commonly interpreted 
as the \textit{physically possible macroscopic states} of very large pieces of matter. 
The extremal decomposition measure $\rho_{\mu}(d \nu)$
for a possibly non-pure state $\mu$ is commonly  interpreted as an expression of a  subjective statistical uncertainty, 
which of the pure states $\nu$ we have in front of us.  \medskip

What is then a (first order) phase transition? 
\textit{A system described by a family of specification kernels is said to have a (first order) phase transition 
if there is more than one infinite-volume Gibbs measure. } This means that the physical system may 
take more than one macroscopic state, for the same values of the parameter, as it is the case for instance 
for large enough $\beta$ in zero magnetic field in the Ising model for $d\geq 2$. 
\medskip

\textit{Lattice Ising model. } The two-dimensional Ising model in zero external field,  has only two extremal Gibbs measures, a plus-like 
state $\mu^+$ and a minus-like state $\mu^-$. At temperatures strictly below the transition temperature 
these states are different, otherwise they coincide. There are no other states, 
which is the content of the Aizenman-Higuchi theorem \cite{Ai80, GeHi00}.  
The Gibbs simplex is richer  for the Ising model on $\Z^3$, where there are also infinite-volume 
states describing stable interfaces at low temperatures, the so-called Dobrushin states \cite{Do72,BrLePf79}. These two-dimensional 
interfaces are analogues of the one-dimensional original boundary curves (Peierls contours) 
in the two-dimensional Ising model we mentioned at the very beginning.

\textit{Inhomogeneous glassy states, Ising model on trees. } Surprising things can happen for spin models like the Ising model on 
infinite graphs (networks) when the graphs are of a much different nature than lattices.  
The DLR formalism makes perfect sense for any countably infinite graph, including for instance trees.  
Models on infinite trees \cite{LyPe16, Ro13}
are mathematically very interesting, and they are often important local building blocks of possibly more complicated 
random networks \cite{Ho17}. 
Consider now the Gibbs simplex in the DLR formalism for 
the nearest neighbor Ising model on a regular infinite tree where each vertex has 
precisely $k+1$ nearest neighbors, $k\geq 2$.  
Here the infinite-volume Gibbs measure which is obtained with open boundary 
conditions is non-pure at very low temperature, and it 
decomposes over uncountably many pure states, which are spatially inhomogeneous. It can moreover be shown  that the 
decomposition measure is atomless at low temperature, which means that  
$\rho_{\mu^\text{free}}(\{\nu\})=0$ for any of the uncountably many pure states $\nu$ \cite{GaMaRuSh20}. It can also be shown that this is a general phenomenon 
for similar models, see  \cite{CoKuLe24} which used the method 
of branch-overlaps motivated by spin-glasses. 
As a loose heuristic why this could be possible, note that  finite subvolumes $\Lambda$ 
of trees have large outer boundaries (which are of a  size of the order of the volume itself). This 
makes it easier (than it would be on lattices) for boundary conditions $\omega$ to be remembered 
by the finite volume Gibbs measure $\gamma_{\Lambda}( d\sigma| \omega)$, even if one is looking for 
the behavior of the spin-configuration $\sigma$ deep inside $\Lambda$. Correspondingly 
there could be and there really are provably (many) more pure states on trees than there are on lattices. 

The appearance of many spatially inhomogeneous pure states is also 
expected to hold on lattices of large enough spatial dimension, however not for 
ferromagnetic models with homogeneous nearest neighbor couplings, 
but for \textit{spin-glass models}. In spin-glasses   the couplings between pairs of spins are random variables 
with symmetric distribution,  see Section \ref{sec:disorder} of this overview.   \medskip 

{\bf Benefits and further properties of infinite-volume measures.} 
The DLR-infinite volume formalism is the natural setup to study  the following subjects. 
There is usually a  model independent part of the theory \cite{Ge11,FrVe17}
which has to be supplemented by model-dependent studies for concrete spin 
systems like the Ising model and its relatives. 

1. \textit{Symmetries of Gibbs measures} either under spin-space transformations or 
geometric transformations.  Think e.g. of lattice translations 
on $\Z^d$  or more generally of a graph automorphism of the underlying graph, which acts on spin configurations, 
and correspondingly on measures on spin configurations. Clearly any boundary condition in the two-dimensional 
Ising model breaks the translation symmetry, but translational symmetry
is restored in the infinite volume, e.g. for the infinite-volume state $\mu^+$. This is relatively easy to show using 
arguments of stochastic monotonicity which hold e.g. for the Ising model \cite{GeHeMa01,FrVe17}. 
Correspondingly one may study ergodicity properties of the 
shift-invariant Gibbs measures, which are invariant  under the translation group on $\Z^d$. \medskip 

2. \textit{Large deviations, concentration properties of observables, Gibbs variational principle in infinite volume}. This can be studied both 
 in the simpler uniqueness regime, and in phase transition regimes. 
 \textit{Large deviations theory} studies  i.e. asymptotics for $\mu$-probabilities of deviations from the expected value w.r.t.  
 a Gibbs measure $\mu$
of an observable like 
$\frac{1}{|\Lambda_n|}\sum_{i \in \Lambda_n}\sigma_i$, but it is much more general.  
 A pioneering figure in the study of large deviations, in particular 
for spatial models, is 
S.R.S. Varadhan \cite{Va84} (born 1940 Madras, India, 
Abel prize 2007). Large deviations theory is intimately connected 
with thermodynamic potentials (like the free energy) and their Legendre transforms, and in this way, it is tied to classical concepts 
in thermodynamics \cite{DeZe10,El06}. 
Legendre transforms of free energies reappear as 
\textit{large deviation rate functionals} (which quantify the exponential decay rate of improbable events).   
Properties depend strongly 
whether we are in the uniqueness regime at small temperature, 
or in the non-uniqueness regime at low temperatures.  

The study of \textit{concentration properties} of functionals of dependent random variables (starting 
from independence) is 
an important subject in many subfields of probability 
\cite{Le01, Ta96}. 
It is related and more general in its character,  
but different from the desire to study precise large deviations. 
In the example of the Ising model one 
is interested in the control of $\mu$-probabilities of deviations of 
the function $F((\sigma_{i})_{\in W})$ depending on spins in the volume $W$, from its expectation,  
w.r.t. a Gibbs measure $\mu$. One expects to find bounds which are exponentially small in $|W|$, depending 
on the form of $F$, at least in uniqueness regimes, 
and weaker bounds in non-uniqueness regimes \cite{ChCoKuRe07}. 
Opposed to large deviations, no symmetries are assumed for $F$ and only quantities measuring 
the overall 
influence of the individual random variables  appear. 

\textit{Large deviations under a Gibbs measure} $\mu(d\omega)$ like the lattice Ising model at any 
temperature 
can be abstractly studied with a high level perspective, namely on the level of \textit{spatial empirical measures}  
in large boxes $\Lambda$. These empirical measures 
are infinite-volume measures themselves, which are of the form $\frac{1}{|\Lambda|}\sum_{i\in \Lambda}\delta_{\tau_{-i} \omega}$, 
where $\tau_{i}$ is a lattice-shift. Here the large deviation functionals appearing which quantify 
the atypicality of a given infinite-volume lattice 
measure, 
turn out to be intimately connected with the Gibbs variational principle in infinite-volume. They 
are connected via the concept of relative 
entropy per site, see 
\cite{Ge11,RaSe15, El06}. 
The \textit{Variational principle in infinite volume} identifies the (possibly multiple!) 
translation-invariant Gibbs measures for a given specification as those having 
the smallest free energy (per lattice site) relative to the interaction.  
This assumes sufficient regularity for the interactions. 
Many of these general theories (for instance the variational principle) assume a certain regularity (which for an Ising model means: enough locality) of the interactions between the spin variables $\sigma_i$ and $\sigma_j$ of the system, as cornerstone theorems 
of the theory may else break down, as examples show \cite{KuLeRe04}. 
The renormalization group pathologies sketched in Section 4 lead to the failure of this necessary regularity in many cases, 
and therefore so-called non-Gibbsian measures with unusual and interesting properties may appear which behave differently then standard short-range theory predicts.

\section{Scaling limits of critical systems in d = 2, interfaces, 
and stochastic Loewner evolution} 
\label{sec:SLE}

Stochastic Loewner evolution (SLE) was invented in 2000 
by Oded Schramm (1961 Jerusalem - 2008 Washington) \cite{Sc00}.  It is a one-parameter class of distributions on fractal planar curves in a given complex domain, 
e.g. in the upper half of the complex plane \cite{Ke17}. It is indexed by one real parameter 
$\kappa$ governing its properties. 
Schramm's brilliant idea was to reconsider the deterministic Loewner differential 
equation of complex analysis \cite{Lo23}  (which was discovered independently at the time  
while Ising was working on his Ph.D. thesis).  It is used to construct conformal maps 
between various domains. Schramm added a Brownian motion, with 
a coupling strength $\kappa$, to map the upper half plane to a domain 
with a random slit, which is equivalently described by a curve. 
He obtained in this way a family of random fractal curves, which are very different from 
diffusion paths, but instead have a certain domain Markov property.  

SLE (which is also called \textit{Schramm} Loewner evolution) is of relevance as describing the possible distributional scaling limits 
for models of two-dimensional random curves at the critical point. 
A major example is the boundary curve of spin-clusters 
for the two-dimensional Ising model in zero external field. A good setup for this on the discrete side is to consider the upper half lattice, where 
it is more convenient to choose the hexagonal lattice then the square lattice, and where the boundary curve 
can be forced into the system as follows. Consider the standard two-dimensional Ising model, together with
spin boundary condition on the $x$-axis which changes from spin value minus (for boundary sites left to the origin) 
to spin-value plus  (for sites right or equal to the origin).  
Then, for each spin configuration  observe the curve which separates plusses from minusses in the bulk, which is anchored at the  
the separation point between plusses and minusses on the boundary line. This curve is random  w.r.t. 
the infinite-volume Gibbs measure on the upper half lattice.  
This Gibbs measure is indeed unique when the inverse temperature $\beta$ equals the critical value $\beta_c$, 
and for this value the scaling limit of the curve has a very interesting self-similar behavior described by SLE \cite{ChDuHoKeSt14}.  
In this case of the exploration curve in the Ising model, and for many more critical two-dimensional models,  the SLE approach is able to obtain 
explicit values of critical exponents in dependence of $\kappa$. These were sometimes promised before by non-rigorous renormalization group 
theory (and previous non-rigorous theories of conformal invariance of theoretical physics), but they are now 
obtained  in a fully rigorous way, which does 
not rely 
on explicit solvability of a model.  
Solvability of a model of course holds 
only for very few models.  

A number of Fields medals were awarded for research on topics 
related to SLE and mathematical statistical mechanics. 
Wendelin Werner (born 1968 in Cologne) obtained the Fields medal in 
2006.  
Stanislav Smirnov (born 1970 in Leningrad) \cite{Sm01} obtained the Fields medal in 
2010. 
Hugo Duminil-Copin (born 1985 in Ile-de-France) \cite{DuTa16}
obtained the Fields medal 2022 for related work, 
but with more focus also on dimensions $d = 3, 4$ where 
SLE is not available. 

What about higher dimensions? We note that in more than two dimensions, interfacial 
\textit{lines} in the Ising model become \textit{interfaces}. 
Such interfaces have 
an effective probabilistic 
description in terms of \textit{gradient models}, whose 
study is a beautiful and rich branch of statistical mechanics. 
In gradient models 
Ising spins are replaced by real-valued 
or integer-valued \textit{height-variables} $\phi_i$ describing 
the height of an interface at a site $i$ in a reference base-plane. 
The interaction 
depends only on the collection of increments $\phi_i-\phi_j$
instead of values of total heights. 
For the particular case of quadratic interactions the distribution becomes 
the \textit{Gaussian free field} (on the lattice), which has a continuum 
limit which also appears as CLT-type scaling limit of a class of non-Gaussian gradient models. Gradient models pose challenging 
questions around the topics of localization and delocalization
\cite{Ve06}, scaling 
limits \cite{NaSp97}, (in-)stability against stochastic perturbations describing quenched disorder 
(in the sense of Section \ref{sec:disorder}) \cite{DaMaPe23, CoKu15}. They are also studied on 
more general graphs than lattices as a base-plane where new phenomena of 
phase coexistence appear \cite{HeKu21}. 

Coming back to two lattice dimensions, SLE can be upgraded to describe not only random planar curves 
anchored at a point, but whole families of random loops, the so-called 
\textit{conformal loop ensembles} describing loop soups. 
It is further used e.g. in the context of 
the Gaussian free field \cite{WePo21}, 
and models of quantum gravity \cite{MiShWe22} where one 
considers distributions of random graphs and their continuum limits which are physically motivated as a description of  fluctuating geometries of space times.


\section{Disordered Ising models} 
\label{sec:disorder}
Consider a 
random \textit{quenched disordered} Gibbs distribution 
for Ising spins $\sigma_i$ taking values $\pm 1$, 
with spin probabilities in finite volume $\Lambda$ according 
to $$\frac{1}{Z_{\Lambda}}\exp(\beta \sum_{i \sim j}J_{ij} \sigma_i \sigma_j +\sum_{i} h_i \sigma_i ).$$ 
Two types of random interactions may appear. 
In \textit{quenched random field} systems the local magnetic fields $h_i$
are chosen according to some a priori distribution $\P$, 
often i.i.d., and then kept fixed. 
\textit{Quenched} physically means that the local degrees of freedom creating disorder 
move on a much larger time-scale and appear to be frozen for the much 
more rapidly fluctuating magnetic moments, and also for the human observer. 
For this reason the further analysis of the structure of Gibbs measures and their 
properties has to be performed conditional on typical realizations of 
the disorder variables. Mathematically speaking, this conditional procedure is the \textit{definition} of a quenched model
\cite{Bo06}.
While translation-invariance is $\P$-a.s. destroyed for the distribution of the spin-variables 
\text{translation-covariance} often holds and ergodicity arguments w.r.t. 
to the disorder distribution are useful. 
Translation covariance of a measure means invariance under joint lattice shift of 
disorder variables and spin variables. In other words, the functional dependence of a Gibbs measure 
on variations of (for instance) the collections of random fields 
$(h_i)_{i\in \Z^d}$ is the same for all regions in space.  

Similarly, the collection of random interactions $J_{ij}$ may be chosen 
according to a quenched distribution $\P$, often independently over 
the edges of the graph under consideration. 
Spin glasses, 
as the most challenging of the classes of random models 
discussed here, are defined by the property that 
positive and negative couplings appear with the same 
probability. 
This is in particular  the case 
for the nearest neighbor Gaussian spinglass model \cite{EdAn75} proposed in the 1970s by
Edwards and Anderson (Nobel prize in physics 1977) 
where the couplings $J_{ij}\sim {\cal{N}}(0,1)$ are standard normal random variables, independently over pairs of nearest neighbors $i,j$ on the lattice. 
Motivation for such models and their relatives is derived from physics and material science. 
There is much interest also from the fields of \textit{network science} \cite{Ho17}, 
and for related models from the fields of 
\textit{computer science and algorithms}, and \textit{machine learning}, see \cite{MeMo09}.

The questions to be asked in the framework 
of infinite volumes Gibbs measures, like the Ising model 
on a $d$-dimensional lattice, are the following. 

{\bf (i)} When 
does randomness of the interaction \textbf{destroy a phase transition}?  

{\bf (ii)} When 
does randomness \textbf{alter an existing phase transition} but leaves it mainly intact?  

{\bf (iii)} When does randomness of the interaction \textbf{create new complex phases}?  

For case {\bf (i)} the random field Ising model with 
deterministic couplings $J$ 
in spatial dimensions $d=2$ is a famous example. Here the  
randomness of the local magnetic field 
terms destroys the low temperature Ising phase transition. 
Following physical predictions of Imry and Ma \cite{ImMa75},  
Aizenman and Wehr 
\cite{AiWe90}
proved that the Gibbs measure at any temperature, 
for arbitrarily weak i.i.d. random fields 
must be almost surely unique, hence the phase transition in the pure model is destroyed. 
They used in their proof the ergodicity of the 
random field distribution, together with a surface vs. volume 
argument which is derived from energetic considerations for  
Peierls contours in the presence of external fields, and is already the core 
of the heuristics of Imry and Ma \cite{ImMa75}. 
In their paper \cite{AiWe90} they also introduced the fruitful notion 
of a \textit{metastate} $\kappa^{\eta}(d\mu)$. A metastate is a random measure on infinite-volume measures. It is random as it depends 
on the given realization of the disorder variables $\eta$, 
here given by the collection of random fields. 
For fixed $\eta$ it is then a probability measure giving weights to the possible
Gibbs measures $\mu$ in the sense of the DLR theory, which lie in the  $\eta$-dependent 
Gibbs simplex. 
This notion has been extended to describe and quantify \textit{chaotic side dependence} 
by Newman and Stein, see \cite{NeSt13,NeReSt22, Ku97,Bo06, Ko23, }.  Chaotic size dependence means that 
a sequence of finite-volume Gibbs measures $\mu^\eta_{\Lambda_n}$, taken 
along an increasing sequence of cubes $\Lambda_n\uparrow \Z^d$ for fixed disorder realization $\eta$, oscillates chaotically between various Gibbs measures, while 
it approaches the corresponding Gibbs simplex. The Newman-Stein metastate is then the large-volume 
limit of the empirical distribution of Gibbs measures at fixed $\eta$, along a volume-sequence. 
Hence it is a measure of relevance of the possible Gibbs-measures, when a large but finite system is considered. 

Proofs of the Aizenman-Wehr phenomenon of washing out a phase transition 
have been given for  
a number of systems, including systems with continuous symmetries. 
The general idea is that quenched 
randomness makes it harder for a system to build up long-range order, so higher dimensions are needed to find a non-uniqueness regime for Gibbs measures. 
There is recent progress in \textit{quantifications} of this 
phenomenon for the two dimensional random field Ising model \cite{AiHaPe20,DiLiXi24}. 
They show that the \textit{influence of choosing different boundary conditions} to the 
expectation of the spin variable at the origin, taken in the corresponding finite-volume Gibbs measures in a box of sidelength $L$, 
decays exponentially fast in $L$.  

The Ising model in small quenched random fields provides also an example for 
case (ii) as it preserves the phase transition at low temperatures in $d\geq 3$. Indeed, it 
does have a ferromagnetic phase transition 
at low temperatures in small disorder, which means for instance 
centered Gaussian random fields with small enough variance. This was under dispute 
in the theoretical physics community, but the issue was famously resolved  
by the renormalization group proof of Bricmont and Kupiainen \cite{BrKu88} on the level of contours. This method is 
conceptually beautiful but technically demanding. A simpler and shorter 
proof appeared much later in \cite{DiZh24}. 

The fundamental questions after 
stability of the phase transition have analogues for models on large but finite graphs. 
An extreme and supposedly simple case is  
the \textit{complete graph} where all pairs of sites are connected. 
Such models on the complete graph are called \textit{mean-field models}, 
since each spin variable experiences an effective field (mean field) created to equal 
parts by all of the other spins, 
which should lead to selfconsistency equations helping the analysis.  
This can be made rigorous relatively easily 
for the non-random mean field Ising model (also known as Curie-Weiss model) \cite{FrVe17}. 
The Curie-Weiss model is a textbook example which is perfectly understood 
to fine probabilistic details including for instance non-trivial scaling laws 
at the critical points.  The analogous questions for 
random mean field models, in particular those of spinglass type, are very hard. 
Let us consider the mean-field spin-glass introduced 
by Sherrington and Kirkpatrick \cite{ShKi75}. They
introduced it already in the title of their paper as a 
\textit{solvable model} but without giving a solution.  
In this model all pairs of $N$ Ising spin variables, interact via a coupling 
$J_{ij}$ which is a centered Gaussian with variance $\frac{1}{N}$, 
chosen independently over the pairs. 

In a branch closely related to spin glasses, 
coupled Ising spins can also serve as an artifical neural network, when they are coupled 
with specifically chosen prescribed constants $J_{ij}$. These can be chosen as to allow "to memorize" 
previously given patterns, in the sense that the Gibbs distribution gives high weights 
to spin configurations which resemble these patterns. 
This was the idea of the physicist John 
Hopfield (born 1933 in Chicago) \cite{Ho82}, 
who received a Nobel prize in physics for his foundational 
contributions to machine learning in 2024, and started a fruitful line of research. 

Giorgio Parisi (born 1948 in Rome) gave a 
(mathematically non-rigorous) solution of the quenched mean-field spin-glass, in the sense 
that he provided the so-called Parisi-formula for the free energy. 
He was awarded the  
 Nobel Prize in Physics (2021) for this work and more generally his work on the "interplay of disorder and fluctuations in physical systems". 
 In his Nobel prize lecture (available on youtube)
 Parisi charmingly recalls the time and situation he found the solution by means of (in his own words) "crazy and illegal" 
 operations, involving the Replica-trick and replica symmetry breaking.  
 
A rigorous breakthrough in the understanding of mean-field spin glasses 
was achieved by Michel Talagrand (born 1952 in B\'eziers, France). 
He was able to prove the validity of Parisi Formula which was published in 2006 \cite{Ta06}. 
See also the complementary 
\cite{Pa12, Pa13} for insights on the level of the random Gibbs measures itself. 
Talagrand is also known in probability among other things for 
his work about the concentration of measure phenomenon, and chaining. The latter 
    is a method to control suprema of dependent processes. 
Talagrand was awarded the Abel Prize in 2024, for his contributions to "probability theory, 
functional analysis, and outstanding applications to mathematical physics and statistics".

A world of related questions appears when we realize that 
quenched randomness may also come from an underlying random geometry.       
What can we say then about the influence of such types of 
\textit{spatial disorder} on collective spin-behavior 
(think of dilute lattices, random networks, point clouds in the continuum, 
dynamical networks, ... ) ? 
In this review we restricted ourselves to static properties, but there is even another 
world of related questions when we begin to talk about associated spin dynamics. 
\section*{Declarations}

\textbf{Data Availability:} Data sharing is not applicable to this article. \\
\textbf{Competing interests:} The author 
has no competing interests to declare that are relevant to the content of this article.
\sloppy
\printbibliography

\end{document}